\documentclass[twoside,12pt]{article}
\usepackage[]{amsmath,amsthm,amssymb}
\usepackage[latin1]{inputenc}

\begin{document}

\title{{\Large\bf  Discrete index Whittaker  transforms}}

\author{Semyon  YAKUBOVICH}
\maketitle

\markboth{\rm \centerline{ Semyon   YAKUBOVICH}}{}
\markright{\rm \centerline{INDEX WHITTAKER  TRANSFORM}}

\begin{abstract} {\noindent Discrete analogs of the index Whittaker  transform are introduced and investigated. It involves series and integrals with respect to a second parameter of the Whittaker  function $W_{\mu, {i n} }(x), \ x >0, \  \mu \in \mathbb{R}, \ n \in \mathbb{N}, \ i $ is the imaginary unit. The corresponding inversion formulas for suitable  functions and sequences in terms of these series and integrals are established.  }

\end{abstract}
\vspace{4mm}



\vspace{4mm}

\section {Introduction and preliminary results}

In 1964 Wimp [3] discovered a reciprocal pair of integral transformations, involving the Whittaker function $ W_{\mu,i\tau}(x)$ [2], Vol. III 

$$F(\tau)=  \int_0^\infty   W_{\mu,i\tau}(x)  f(x) {dx\over x^2},\  \tau >0,\eqno(1.1)$$

$$f(x)= { 1 \over \pi^2  } \int_0^\infty  \tau \sinh(2\pi\tau) \left|\Gamma\left({1\over 2}- \mu + i\tau\right)\right|^2  W_{\mu,i\tau}(x)  F(\tau) d\tau,\ x >0,\eqno(1.2)$$
 where $\mu \in \mathbb{R}$ and $\Gamma(z)$ is the Euler gamma function [2], Vol. III.  These transforms are particular cases of a more general integral transform with Meijer $G$-function as the kernel [4].  As we see in (1.1), (1.2), the corresponding integrals depend upon a second parameter (index) of the Whittaker function and involve it as the variable of integration. Therefore formulas (1.1), (1.2) are  called the reciprocal index Whittaker transforms.  When $\mu=0$ this pair can be reduced to the familiar Kontorovich-Lebedev transforms [6]. 
 
 In this  paper we will deal with the following discrete analogs of the index Whittaker transforms
 
 $$f(x)=  e^{-x/2} \sum_{n=1}^\infty  a_n\  W_{\mu, {i n}}(x),\quad x >0,\eqno(1.3)$$
 
 $$a_n =   \int_0^\infty   e^{-x/2} W_{\mu, {i n\over 2}}(x)  f(x) x^{\mu-2} dx,\  n \in \mathbb{N}\eqno(1.4)$$
 for  suitable classes of functions $f$ and sequences $\{a_n\}_{n\ge 1}$.   Our approach will be based on the use of classical Fourier series and some integrals, involving the Whittaker and modified Bessel functions [2], Vol. II.  As is known,  the Whittaker function has the following asymptotic behavior at the origin and near infinity
 
 $$W_{\mu,\nu}(x)= O\left( x^{1/2 + {\rm Re} \nu}\right) +   O\left( x^{1/2 - {\rm Re} \nu}\right),\quad  x \to 0+,\eqno(1.5)$$
 
 $$ W_{\mu,\nu}(x)= O\left( e^{-x/2} x^{{\rm Re} \mu}\right),\ x \to \infty.\eqno(1.6)$$
 It can be represented by the  integral [5]

 $$ \Gamma\left({1\over 2}- \mu + i\tau\right)\Gamma\left({1\over 2}- \mu - i\tau\right) W_{\mu,i\tau}(x) $$
 
 $$= 2(4x)^\mu e^{-x/2} \int_0^\infty  t^{-2\mu} e^{-t^2/(4x)} K_{2i\tau}(t) dt,\ x >0, \  {\rm Re} \mu < {1\over 2},\ \tau \in \mathbb{R},\eqno(1.7)$$
 where $K_\nu(z)$ is the modified Bessel function.   For the pure imaginary index this function satisfies the inequality (see [6], formula (1.100))
 
  $$\left|K_{i\tau}(x)\right|  \le e^{-\delta|\tau|} K_0\left(x\cos(\delta)\right),\ x >0, \ \tau \in \mathbb{R}, \ \delta \in \left[ 0, {\pi\over 2} \right).\eqno(1.8)$$
Hence, as a consequence of (1.7), (1.8)  the following inequality for the Whittaker function holds 
 
 $$ \left| W_{\mu,i\tau}(x)\right| \le  \left(\cos(\delta)\right)^{-1} \left| {\Gamma \left({1\over 2} - \mu\right)\over  \Gamma\left({1\over 2}- \mu + i\tau\right)}\right|^2 W_{\mu,0}\left(x\cos^2(\delta)\right)
 e^{-{x\sin^2(\delta)\over 2} - 2\delta |\tau|},\eqno(1.9)$$
 where $x >0, \ \tau \in \mathbb{R},\  \mu < {1\over 2}, \delta \in \left[0, {\pi\over 2}\right)$.  On the other hand,  the Whittaker function possesses  by  the Mellin-Barnes representation   in terms of the ratio of  products of gamma functions (see [2], Vol. III, Entry 8.4.44.2)
 
 $$e^{-x/2} W_{\mu,\rho}\left( x\right) $$
 
 $$= {1\over  2\pi i} \int_{\gamma-i\infty}^{\gamma +i\infty} \frac{\Gamma\left({1\over 2} +\rho + s\right) \Gamma\left({1\over 2} -\rho + s\right)}
 {\Gamma(1-\mu+s)} x^{-s} ds,\quad x >0, \ \gamma >   \left|  {\rm Re} \rho \right| - {1\over 2}.\eqno(1.10)$$
 Hence, taking  integral (1.10), we deduce

  $$  \int_0^\infty  e^{-x^2/(4t)- t/2}\  W_{\mu,\rho}(t)\   t^{\mu-2}  dt =  {1\over  2\pi i} \int_{\gamma-i\infty}^{\gamma +i\infty} \frac{\Gamma\left({1\over 2} +\rho + s\right) \Gamma\left({1\over 2} -\rho + s\right)} {\Gamma(1-\mu+s)} $$
  
  $$\times  \int_0^\infty  e^{-x^2/(4t)} t^{\mu-s-2} dt ds,\eqno(1.11)$$
 where the interchange of the order of integration is allowed by Fubini's theorem owing to the estimate
  
$$ \int_{\gamma-i\infty}^{\gamma +i\infty} \left| \frac{\Gamma\left({1\over 2} +\rho + s\right) \Gamma\left({1\over 2} -\rho + s\right)} {\Gamma(1-\mu+s)} \right| \int_0^\infty  e^{-x^2/(4t)} \left|t^{\mu-s-2} dt ds\right|  $$

 $$= \left({x\over 2}\right)^{2({\rm Re} \mu-\gamma-1)} \Gamma(\gamma-  {\rm Re} \mu+1)  \int_{\gamma-i\infty}^{\gamma +i\infty} \left| \frac{\Gamma\left({1\over 2} +\rho + s\right) \Gamma\left({1\over 2} -\rho + s\right)} {\Gamma(1-\mu+s)}  ds \right| < \infty,$$ 
 where $x >0, \gamma >  {\rm Re} \mu-1,\   \left|  {\rm Re} \rho \right| - {1\over 2}$,  and the Stirling asymptotic formula for the gamma function (see  [4], formula  (1.12)). Consequently, returning to (1.11) and calculating the Euler integral, we derive
 
 $$   \int_0^\infty  e^{-x^2/(4t)- t/2}\  W_{\mu,\rho}(t)\   t^{\mu-2}  dt =   \left({x\over 2}\right)^{2( \mu-1)}  {1\over  2\pi i} \int_{\gamma-i\infty}^{\gamma +i\infty} \Gamma\left({1\over 2} +\rho + s\right) $$
 
 $$\times \Gamma\left({1\over 2} -\rho + s\right)  \left({x\over 2}\right)^{- 2s}  ds.$$
 Now, employing the formula (cf. [2], Vol. III, Entry 8.4.23.1)
 
 $$ {1\over  2\pi i} \int_{\gamma-i\infty}^{\gamma +i\infty} \Gamma\left(s +\rho\right) \Gamma\left(s -\rho \right)  x^{-s} ds = 2 K_{2\rho}(2\sqrt x),\ x >0,  \gamma >  \left|  {\rm Re} \rho \right|,$$
 we easily end up with the equality (cf. [2], Vol. III, Entry 2.19.4.7)
 
 $$  \int_0^\infty  e^{-x^2/(4t)- t/2}\  W_{\mu,\rho}(t)\   t^{\mu-2}  dt = 2  \left({x\over 2}\right)^{2\mu-1} K_{2\rho} (x),\ x >0, \mu, \rho \in \mathbb{C},\eqno(1.12)$$
which will be used in the sequel.  Finally in this section we mention a useful integral (cf. [2], Vol. I, Entry 2.3.15.3) in our investigation which defines the parabolic cylinder function $D_\nu(z)$

$$\int_0^\infty x^{\alpha-1} e^{- p x^2-qx} dx = {\Gamma(\alpha) \over (2p)^{\alpha/2}}\  e^{{q^2\over 8p}}\  D_{-\alpha} \left({q\over \sqrt{2p}}\right), \ {\rm Re} \alpha,  \ {\rm Re} p > 0, \ q \in \mathbb{C}.\eqno(1.13)$$
 These preliminary results are, indeed,  key ingredients  to achieve our goal to invert discrete index Whittaker transforms (1.3), (1.4).  The suggested approach is based on reducing to the discrete Kontorovich-Lebedev transform whose theory is recently elaborated by the author in [7].

\section{Inversion theorems}

We begin with

{\bf Theorem 1.} {\it Let   $ \mu < {1\over 2}$ and  the sequence $\{a_n\}_{n\ge 1}$ satisfy  the condition 

$$ \sum_{m=1}^\infty {\left| a_m \right| e^{- 2\delta m} \over \left|\Gamma\left( 1/2 +im - \mu\right)\right|^2 }   < \infty,\quad  \delta \in \left[0, \  {\pi\over 2} \right).\eqno(2.1)$$
Then the discrete transformation $(1.3)$ can be inverted by the formula

$$a_n=    {  2^{1/2+\mu}  \over  \pi^2}\  \Gamma(1-2\mu)   n \sinh(2\pi n) \int_0^\infty  \Phi^\mu_n(t) f(t) t^{-3/2}  dt,\eqno(2.2)$$
where

$$  \Phi^\mu_n(x) =  \int_{0}^\pi  e^{x\cosh^2(u)/2} D_{2\mu-1}\left(\sqrt{2x} \cosh(u)\right) \cos(2nu) du\eqno(2.3)$$
and integral $(2.2)$ converges absolutely.}

\begin{proof}   Taking the modified Laplace transform [4] of both sides of (1.3), we interchange the order of integration and summation and then appeal to (1.12) to derive

$$  \int_0^\infty  e^{-x^2/(4t)} f(t) t^{\mu-2}  dt = \sum_{m=1}^\infty  a_m\   \int_0^\infty  e^{-x^2/(4t)- t/2}   W_{\mu, i m}(t) t^{\mu-2}  dt $$

$$ =   2  \left({x\over 2}\right)^{2\mu-1}  \sum_{m=1}^\infty  a_m  K_{2im} (x),\quad x >0.\eqno(2.4)$$
This interchange is permitted due to Fubini's theorem via inequality (1.9), representation (1.12)  and condition (2.1). In fact, we have

$$  \int_0^\infty  e^{-x^2/(4t)- t/2}  t^{\mu-2} \sum_{m=1}^\infty  \left| a_m\  W_{\mu, i m}(t)\right|  dt \le  \left(\cos(\delta)\right)^{-1} \left[ \Gamma \left({1\over 2} - \mu\right)\right]^2 $$

$$\times  \int_0^\infty  e^{-x^2/(4t)- t(1+\sin^2(\delta))/2}   W_{\mu,0}\left(t\cos^2(\delta)\right)  t^{\mu-2} dt \sum_{m=1}^\infty  {\left| a_m \right|  e^{ - 2\delta m} \over  \left|\Gamma\left(1/2+im - \mu\right)\right|^2 }$$

$$ \le   \left(\cos(\delta)\right)^{-1} \left[ \Gamma \left({1\over 2} - \mu\right)\right]^2 $$

$$\times  \int_0^\infty  e^{-x^2/(4t)- t(\cos^2(\delta))/2}   W_{\mu,0}\left(t\cos^2(\delta)\right)  t^{\mu-2} dt \sum_{m=1}^\infty  {\left| a_m \right|  e^{ - 2\delta m} \over  \left|\Gamma\left(1/2+im - \mu\right)\right|^2 }$$

$$= 2  \left({x\over 2}\right)^{2\mu-1}  \left[ \Gamma \left({1\over 2} - \mu\right)\right]^2 K_0\left(x\cos(\delta)\right)  \sum_{m=1}^\infty  {\left| a_m \right|  e^{ - 2\delta m} \over  \left|\Gamma\left(1/2+im - \mu\right)\right|^2 }  < \infty.$$
Now, the right-hand side of the latter equality in (2.4) is a discrete Kontorovich-Lebedev transform of double index (cf. [7], Th. 1),  which  can be inverted under the condition  $ \sum_{m=1}^\infty | a_m| e^{-\pi m} < \infty.$ But this condition holds by virtue of (2.1) and Stirling's asymptotic formula for the gamma function since $\left|\Gamma\left(1/2+im - \mu\right)\right| = O\left( m^{-\mu} e^{-\pi m/2}\right),\ m \to \infty.$
Consequently, inverting the Kontorovich-Lebedev transform in (2.4), we get, reciprocally,

$$a_n = {4^{\mu} \over \pi^2}\  n \sinh(2\pi n) \int_0^\infty x^{-2\mu} J(x, 2in,\pi)  \int_0^\infty  e^{-x^2/(4t)} f(t) t^{\mu-2}  dt dx,\eqno(2.5)$$
where $J(x, in,\pi)$ is the incomplete modified Bessel function [1]

$$J(x, in,\pi) = \int_0^\pi e^{-x\cosh(u)} \cos(nu) du$$

$$= {x\over n} \int_0^\pi e^{-x\cosh(u)}  \sinh(u) \sin(nu) du,\eqno(2.6)$$
where the second integral in (2.6) is obtained via integration by parts. Hence,  since due to  (1.9) and Entries 2.3.16.1 in [2], Vol. I, 2.16.6.2 in [2], Vol. II we find 

$$ \int_0^\infty x^{-2\mu} \left|J(x, 2in,\pi)\right|   \int_0^\infty  e^{-x^2/(4t)} |f(t)| t^{\mu-2}  dt dx$$

$$\le  {1\over 2n \cos(\delta)}  \left[\Gamma \left({1\over 2} - \mu\right)\right]^2  \int_0^\pi \sinh(u)  \int_0^\infty   \int_0^\infty  \exp\left( - x\cosh(u) -{x^2\over 4t}- {t \over 2} (1+\sin^2(\delta))\right) $$

$$\times   x^{1-2\mu}  W_{\mu,0}\left(t\cos^2(\delta)\right)  t^{\mu-2}  dt dx du  \sum_{m=1}^\infty {\left| a_m \right| e^{- 2\delta m} \over \left|\Gamma\left( 1/2+im - \mu\right)\right|^2 }$$

$$\le  2^{2(1-\mu)}  \left[\Gamma \left({1\over 2} - \mu\right)\right]^2  \int_0^\pi \sinh(u)  \int_0^\infty  e^{-x \cosh(u)} K_{0} \left(x\cos(\delta) \right)  dx du $$

$$\times \sum_{m=1}^\infty {\left| a_m \right| e^{- 2\delta m} \over \left|\Gamma\left( 1/2+im- \mu\right)\right|^2 }$$

$$=   2^{2(1-\mu)}  \left[\Gamma \left({1\over 2} - \mu\right)\right]^2  \int_0^\pi {\sinh(u) \over (\cosh^2(u)-\cos^2(\delta))^{1/2}} $$

$$\times \log\left( {\cosh(u) +  (\cosh^2(u)-\cos^2(\delta))^{1/2}\over \cos(\delta)} \right)du \sum_{m=1}^\infty {\left| a_m \right| e^{- 2\delta m} \over \left|\Gamma\left( 1/2+im - \mu\right)\right|^2 }$$

$$\le    2^{2(1-\mu)} \pi \left[\Gamma \left({1\over 2} - \mu\right)\right]^2  \log\left( {\cosh(\pi) +  (\cosh^2(\pi)-\cos^2(\delta))^{1/2}\over \cos(\delta)}\right) $$

$$\times \sum_{m=1}^\infty {\left| a_m \right| e^{- 2\delta m} \over \left|\Gamma\left( 1/2+im - \mu\right)\right|^2 } < \infty,\eqno(2.7)$$
we interchange  the order of integration in (2.5)  by Fubini's theorem.  Then, appealing to (1.13), (2.6), we arrive at the inversion formula (2.2) with the kernel (2.3), completing the proof of Theorem 1.

\end{proof}

{\bf Remark 1}.  The discrete transformation (1.3) and its inversion formula (2.2) generate the following expansion of an arbitrary sequence $\{a_n\}_{n\ge 1}$, satisfying condition (2.1)

$$a_n=   { 2^{1/2+\mu}  \over  \pi^2}\  \Gamma(1-2\mu)\ n \sinh(2\pi n) \int_0^\infty  \Phi^\mu_n(t)  t^{-3/2}  e^{-t/2} \sum_{m=1}^\infty  a_m W_{\mu, i m}(t) dt.$$

When $\mu=0$, it gives $W_{0, im}(x) = \sqrt{{x\over \pi} } K_{im}\left({x\over 2}\right)$ and (1.3) reduces to the modified discrete Kontorovich-Lebedev transform

$$f(x)=  e^{-x/2}  \sqrt{{x\over \pi} } \sum_{m=1}^\infty  a_m  K_{im}\left({x\over 2}\right),\quad x >0.\eqno(2.8)$$
 Therefore, using the value of the integral (see [2], Vol. I, Entry 2.3.15.4)
 
 $$\int_0^\infty  e^{x^2/(4t)- x\cosh(u)} dx = \sqrt {\pi t}\  e^{t \cosh^2(u)} \hbox{erfc} \left( \sqrt t \cosh(u) \right),\eqno(2.9)$$
where $\hbox{erfc}(z)$ is the complementary error function [2], Vol. I,  we write its inversion formula in the form 

$$a_n=   { n \sinh(2\pi n) \over  \pi\sqrt\pi}  \int_0^\infty  \Phi^0_n(t)  f(t)  t^{-3/2} dt,\quad n \in \mathbb{N},\eqno(2.10)$$
where

$$  \Phi^0_n(x) =  \int_{0}^\pi  e^{x\cosh^2(u)}  \hbox{erfc} \left( \sqrt x \cosh(u) \right) \cos(2nu) du.\eqno(2.11)$$
Concerning the discrete index Whittaker transform (1.4), we have the following result.

{\bf Theorem 2}.   {\it Let  $\mu < 1/2$ and $f$ be a complex-valued function on $\mathbb{R}_+$ which is represented by the integral 

$$f(x) = \Gamma(2(1-\mu))  (2x)^{1-\mu}  \int_{-\pi}^\pi   e^{x\cosh^2(u)/2} D_{2(\mu-1)}\left(\sqrt{2x} \cosh(u)\right)   \varphi(u) du,\quad x >0,\eqno(2.12)$$ 
where $ \varphi(u) = \psi(u)\sinh(u)$ and $\psi$ is a  $2\pi$-periodic function, satisfying the Lipschitz condition on $[-\pi, \pi]$, i.e.

$$\left| \psi(u) - \psi(v)\right| \le C |u-v|, \quad  \forall \  u, v \in  [-\pi, \pi],\eqno(2.13)$$
where $C >0$ is an absolute constant.  Then the following inversion formula for  transformation $(1.4)$  holds

$$ f(x)  =  {  (x/2)^{1-\mu}\over  \pi^2} \  \Gamma(2(1-\mu))   \sum_{n=1}^\infty    \sinh(\pi n)  \Psi^\mu_n (x) a_n,\quad x > 0,\eqno(2.14)$$
where $\Psi^\mu_n (x)$ is defined by }

$$\Psi^\mu_n(x)= \int_{-\pi}^\pi  e^{x\cosh^2(u)/2} D_{2(\mu-1)}\left(\sqrt{2x} \cosh(u)\right) \sinh(u) \sin(nu) du.\eqno(2.15)$$

\begin{proof}  In fact,  appealing to (1.13), we write (2.12) as follows

$$f(x) =  \int_{-\pi}^\pi   \varphi(u) \int_0^\infty e^{-t^2/(4x)- t\cosh(u)} t^{1-2\mu} dt du.\eqno(2.16)$$
Plugging the right-hand side of (2.16) in (1.4), we change the order of integration and calculate the integral with  respect to $x$, employing (1.12),  to obtain 

$$ \int_0^\infty   e^{-x/2} W_{\mu, {i n\over 2}}(x)  f(x) x^{\mu-2} dx = \int_{-\pi}^\pi   \varphi(u) \int_0^\infty e^{- t\cosh(u)} t^{1-2\mu} $$

$$\times  \int_0^\infty   e^{-x/2 -t^2/(4x)} W_{\mu, {i n\over 2}}(x)   x^{\mu-2} dxdt du =  4^{1-\mu} \int_{-\pi}^\pi   \varphi(u) \int_0^\infty e^{- t\cosh(u)} K_{in}(t) dtdu.\eqno(2.17)$$
The interchange of the order of integration is guaranteed by Fubini's theorem via inequality (1.9), the continuity of $\varphi$  and the estimate (cf. (2.7))

$$ \int_0^\infty   e^{-x/2} \left|W_{\mu, {i n\over 2}}(x)\right|  x^{\mu-2} \int_{-\pi}^\pi  | \varphi(u) | \int_0^\infty e^{-t^2/(4x)- t\cosh(u)} t^{1-2\mu} dt du dx$$

$$\le  \left(\cos(\delta)\right)^{-1} e^{-\delta n} \left| {\Gamma \left({1\over 2} - \mu\right)\over  \Gamma\left({1+in\over 2}- \mu \right)}\right|^2 \int_0^\infty   e^{-x \cos^2(\delta) /2} W_{\mu,0}\left(x\cos^2(\delta)\right)  x^{\mu-2}$$

$$\times  \int_{-\pi}^\pi  | \varphi(u) | \int_0^\infty e^{-t^2/(4x)- t\cosh(u)} t^{1-2\mu} dt du dx$$ 

$$=  4^{1-\mu}  e^{-\delta n} \left| {\Gamma \left({1\over 2} - \mu\right)\over  \Gamma\left({1+in\over 2}- \mu \right)}\right|^2  \int_{0}^\pi  | \varphi(u) | \int_0^\infty e^{- t\cosh(u)} K_0(t\cos(\delta))  dt du < \infty.$$ 
Hence, returning to (2.17) and calculating the latter integral by $t$ with the aid of Entry 2.16.6.1 in  [2], Vol. II

$$ \int_0^\infty e^{- t\cosh(u)} K_{in}(t) dt = {\pi \sin(nu) \over \sinh(\pi n) \sinh(u)},$$
we get finally, combining with (1.4) and the definition of $\varphi$,

$$  a_n =   {4^{1-\mu} \pi\over \sinh(\pi n)} \int_{-\pi}^\pi   \varphi(u) {\sin(nu) \over \sinh(u)} du=  {4^{1-\mu} \pi\over \sinh(\pi n)} \int_{-\pi}^\pi   \psi(u) \sin(nu) du.\eqno(2.18)$$
Therefore, following the same scheme as in the proof of Theorem 5 in [7],  we substitute the value of $a_n$ by (2.18) and $\Psi^\mu_n (x)$ by (2.15) into the partial sum of the series (2.14). Then,    calculating this sum via the known identity,  we  obtain 

$$S_N(x)=  { (x/2)^{1-\mu}\over  \pi^2} \  \Gamma(2(1-\mu))   \sum_{n=1}^N    \sinh(\pi n)  \Psi^\mu_n (x) a_n$$

$$ = { (2x)^{1-\mu} \over \pi }\  \Gamma(2(1-\mu))  \sum_{n=1}^N  \int_{-\pi}^\pi   e^{x\cosh^2(t)/2} D_{2(\mu-1)}\left(\sqrt{2x} \cosh(t)\right) \sinh(t)\sin( nt) dt $$

$$\times   \int_{-\pi}^\pi  \psi(u) \  \sin( nu) du$$

$$ = { (2x)^{1-\mu} \over 4\pi }  \Gamma(2(1-\mu)) \int_{-\pi}^\pi    e^{x\cosh^2(t)/2} D_{2(\mu-1)}\left(\sqrt{2x} \cosh(t)\right) \sinh(t) $$

$$\times   \int_{-\pi}^\pi  \left[ \psi(u)- \psi(-u) \right]  \  {\sin \left((2N+1) (u-t)/2 \right)\over \sin( (u-t) /2)}  du dt.\eqno(2.19)$$
Since $\psi$ is $2\pi$-periodic, we treat  the latter integral with respect to $u$ as follows 

$$  \int_{-\pi}^{\pi}  \left[ \psi(u)- \psi(-u) \right]  \  {\sin \left((2N+1) (u-t)/2 \right)\over \sin( (u-t) /2)}  du $$

$$=  \int_{ t-\pi}^{t+ \pi}  \left[ \psi(u)- \psi(-u) \right]  \  {\sin \left((2N+1) (u-t)/2 \right)\over \sin( (u-t) /2)}  du $$

$$=  \int_{ -\pi}^{\pi}  \left[ \psi(u+t)- \psi(-u-t) \right]  \  {\sin \left((2N+1) u/2 \right)\over \sin( u /2)}  du. $$
Moreover,

$$ {1\over 2\pi} \int_{ -\pi}^{\pi}  \left[ \psi(u+t)- \psi(-u-t) \right]  \  {\sin \left((2N+1) u/2 \right)\over \sin( u /2)}  du - \left[ \psi(t)- \psi(-t) \right] $$

$$=  {1\over 2\pi} \int_{ -\pi}^{\pi}  \left[ \psi(u+t)- \psi(t) + \psi (-t) - \psi(-u-t) \right]  \  {\sin \left((2N+1) u/2 \right)\over \sin( u /2)}  du.$$
When  $u+t > \pi$ or  $u+t < -\pi$ then we interpret  the value  $\psi(u+t)- \psi(t)$ by  formulas

$$\psi(u+t)- \psi(t) = \psi(u+t-2\pi)- \psi(t - 2\pi),$$ 

$$\psi(u+t)- \psi(t) = \psi(u+t+ 2\pi)- \psi(t +2\pi),$$ 
respectively.  Analogously, the value  $\psi(-u-t)- \psi(-t)$  can be treated.   Then   due to the Lipschitz condition (2.13) we have the uniform estimate
for any $t \in [-\pi,\pi]$

$${\left|  \psi(u+t)- \psi(t) + \psi (-t) - \psi(-u-t) \right| \over | \sin( u /2) |}  \le 2C \left| {u\over \sin( u /2)} \right|.$$
Therefore,  owing to the Riemann-Lebesgue lemma

$$\lim_{N\to \infty } {1\over 2\pi} \int_{ -\pi}^{\pi}  \left[ \psi(u+t)- \psi(-u-t)  - \psi(t) + \psi (-t) \right]  \  {\sin \left((2N+1) u/2 \right)\over \sin( u /2)}  du =  0\eqno(2.20)$$
for all $ t\in [-\pi,\pi].$    Besides, returning to (2.19), we estimate the iterated integral 

$$ \int_{-\pi}^\pi  e^{x\cosh^2(t)/2} D_{2(\mu-1)}\left(\sqrt{2x} \cosh(t)\right) | \sinh(t)  | $$

$$\times \int_{ -\pi}^{\pi} \left| \left[ \psi(u+t)- \psi(-u-t)  - \psi(t) + \psi (-t) \right]  \  {\sin \left((2N+1) u/2 \right)\over \sin( u /2)}  \right| du dt$$

$$\le  4 C \int_{0}^\pi  e^{x\cosh^2(t)/2} D_{2(\mu-1)}\left(\sqrt{2x} \cosh(t)\right)  \sinh(t)  dt $$

$$\times   \int_{ -\pi}^{\pi}   \left| {u\over \sin( u /2)} \right| du < \infty,\ x >0.$$
Consequently, via  the dominated convergence theorem it is possible to pass to the limit when $N \to \infty$ under the  integral sign, and recalling (2.20), we derive

$$  \lim_{N \to \infty}   { (2x)^{1-\mu}\over 4\pi } \Gamma(2(1-\mu)) \int_{-\pi}^\pi e^{x\cosh^2(t)/2} D_{2(\mu-1)}\left(\sqrt{2x} \cosh(t)\right) \sinh(t)  $$

$$\times  \int_{ -\pi}^{\pi}  \left[ \psi(u+t)- \psi(-u-t)  - \psi(t) + \psi (-t) \right] $$

$$\times  \  {\sin \left((2N+1) u/2 \right)\over \sin( u /2)}  du dt $$

$$=  {(2x)^{1-\mu}  \over 4\pi }   \Gamma(2(1-\mu)) \int_{-\pi}^\pi   e^{x\cosh^2(t)/2} D_{2(\mu-1)}\left(\sqrt{2x} \cosh(t)\right) \sinh(t) $$

$$ \times \lim_{N \to \infty}  \int_{ -\pi}^{\pi}  \left[ \psi(u+t)- \psi(-u-t)  - \psi(t) + \psi (-t) \right]  \  {\sin \left((2N+1) u/2 \right)\over \sin( u /2)}  du dt = 0.$$
Hence, combining with (2.19),  we obtain  by virtue of  the definition of $\varphi$ and $f$

$$ \lim_{N \to \infty}  S_N(x) =   {(2x)^{1-\mu} \over 2}   \Gamma(2(1-\mu))  \int_{-\pi}^\pi  e^{x\cosh^2(t)/2} D_{2(\mu-1)}\left(\sqrt{2x} \cosh(t)\right) $$

$$\times \left[  \varphi (t)+ \varphi(-t) \right] dt = f(x),$$
where the integral (2.12) converges since $\varphi \in C[-\pi,\pi]$.  Thus we established  (2.14), completing the proof of Theorem 2.
 
\end{proof} 

{\bf Remark 2}.  Functions (2.3) and (2.15) are related by the formula $\Psi^\mu_n(x) = 2 {\rm Im} \Phi^{\mu-1/2}_{n-i\over 2} (x)$.

\bigskip
\centerline{{\bf Acknowledgments}}
\bigskip

\noindent The work was partially supported by CMUP, which is financed by national funds through FCT (Portugal)  under the project with reference UIDB/00144/2020.

\bigskip
\centerline{{\bf References}}
\bigskip
\baselineskip=12pt
\medskip
\begin{enumerate}

\item[{\bf 1.}\ ]   D.S. Jones,   Incomplete Bessel functions. I,  {\it  Proc. Edinb. Math. Soc.} {\bf 50} ( 2007),  N 1,  173-183.

\item[{\bf 2.}\ ] A.P. Prudnikov, Yu.A. Brychkov and O.I. Marichev, {\it Integrals and Series}. Vol. I: {\it Elementary Functions}, Vol. II: {\it Special Functions}, Gordon and Breach, New York and London, 1986, Vol. III : {\it More special functions},  Gordon and Breach, New York and London,  1990.

\item[{\bf 3.}\ ]  J. Wimp,  A class of integral transforms,  {\it  Proc. Edinb. Math. Soc.} {\bf 14} (1964), N 2,  33-40. 

\item[{\bf 4.}\ ]  S. Yakubovich and Yu. Luchko, The Hypergeometric Approach to Integral Transforms and Convolutions, {\it Kluwer
Academic Publishers, Mathematics and Applications.} Vol.287, 1994.

\item[{\bf 5.}\ ]  H.M.  Srivastava,  Yu. V. Vasil'ev and  S.B.  Yakubovich,  A class of index transforms with Whittaker's function as the kernel. {\it Quart. J. Math. Oxford} {\bf 49} (2) (1998),  375-394. 

\item[{\bf 6.}\ ] S. Yakubovich, {\it Index Transforms}, World Scientific Publishing Company, Singapore, New Jersey, London and
Hong Kong, 1996.

\item[{\bf 7.}\ ]  S. Yakubovich, Discrete Kontorovich-Lebedev transforms, {\it Ramanujan J.}  DOI 10.1007/s\\ 11139-020-00313-7.

\end{enumerate}

\vspace{5mm}

\noindent S.Yakubovich\\
Department of  Mathematics,\\
Faculty of Sciences,\\
University of Porto,\\
Campo Alegre st., 687\\
4169-007 Porto\\
Portugal\\
E-Mail: syakubov@fc.up.pt\\

\end{document}